\def\versiondate{July 22, 1999 and August 20, 1999}

\documentclass{amsart}
\usepackage{amsmath,amssymb,epsf,amscd}

\newif\ifmultifile
\multifilefalse

\ifmultifile
  \def\myendinput{\endinput}
\else
  \def\myendinput{}
\fi

\def\Complexes{{\mathbf C}}
\def\Field{{\mathbf F}}
\def\Identity{\operatorname{Id}}
\def\Integers{{\mathbf Z}}
\def\mysetminus{\smallsetminus}  
\def\Projective{{\mathbf P}}
\def\projection{\operatorname{pr}}
\def\Reals{{\mathbf R}}

\newtheorem{theorem}{Theorem}[section]
\newtheorem{prop}[theorem]{Proposition}
\newtheorem{lemma}[theorem]{Lemma}
\newtheorem{cor}[theorem]{Corollary}

\theoremstyle{definition}
\newtheorem{definition}[theorem]{Definition}
\newtheorem{example}[theorem]{Example}

\theoremstyle{remark}
\newtheorem{remark}[theorem]{Remark}

\numberwithin{equation}{section}

\begin{document}

\title{Blowup and Fixed Points}

\author{C.~W.~Stark}
\address{Department of Mathematics, University of Florida,
PO Box 118105, Gainesville, Florida 32611--8105}
\curraddr{National Science Foundation,
Division of Mathematical Sciences,
4201 Wilson Boulevard, Room 1025,
Arlington, Virginia 22230}
\email{cws@math.ufl.edu, cstark@nsf.gov}
\thanks{The author was supported in part by an NSA grant.}


\date{\versiondate}

\begin{abstract}
Blowing up a point $p$ in a manifold $M$ builds a new manifold
$\widehat M$ in which $p$ is replaced by the projectivization of the
tangent space $T_pM$.  This well--known operation also applies to
fixed points of diffeomorphisms, yielding continuous homomorphisms
between automorphism groups of $M$ and $\widehat M$.  The construction
for maps involves a loss of regularity and is not unique at the lowest
order of differentiability.  Fixed point sets and other aspects of
blownup dynamics at the singular locus are described in terms of
derivative data; $C^0$ data are not sufficient to determine 
much about these issues.

Topological generalizations of the blowup construction
prove to be much less natural than the classical versions, and
no lifting homomorphism for homeomorphism groups can be constructed.
\end{abstract}

\maketitle

\ifmultifile
\input intro
\input blowups
\input model
\input natural
\input topology
\input maps
\input dynamic
\input variant
\bibliographystyle{amsplain}
\bibliography{papers}

\begin{thebibliography}{1}

\bibitem{AbkulutKing91}
S.~Akbulut and H.~King, \emph{Rational structures on $3$-manifolds}, Pacific J.
  Math. \textbf{150} (1991), 201--214.

\bibitem{BenedettiMarin92}
R.~Benedetti and A.~Marin, \emph{D\'echirures de vari\'et\'es de dimension
  trois et la conjecture de {N}ash de rationalit\'e en dimension trois},
  Comment. Math. Helv. \textbf{67} (1992), 514--545.

\bibitem{Harris92}
J.~Harris, \emph{{Algebraic Geometry: A First Course}}, Graduate Texts in
  Math., no. 133, Springer--Verlag, New York, 1992.

\bibitem{KatokHasselblatt95}
A.~Katok and B.~Hasselblatt, \emph{{Introduction to the Modern Theory of
  Dynamical Systems}}, Encyclopedia of Mathematics and its Applications,
  no.~54, Cambridge University Press, 1995.

\bibitem{McDuffSalamon98}
D.~McDuff and D.~Salamon, \emph{{ Introduction to Symplectic Topology, Second
  Edition}}, Oxford University Press, 1998.

\bibitem{Mikhalkin97}
G.~Mikhalkin, \emph{Blowup equivalence of smooth closed manifolds}, Topology
  \textbf{36} (1997), 287--299.

\bibitem{Nash52}
J.~Nash, \emph{Real algebraic manifolds}, Ann. of Math. (2) \textbf{56} (1952),
  405--421.

\end{thebibliography}
\fi

\empty


\section{Introduction}
The construction for blowing up points and subspaces which 
is a mainstay in algebraic geometry, especially in the resolution
of singularities, is investigated here from a dynamical point of
view.  The blowup in this sense of a point $p$ in a smooth manifold
$M^n$ is a map of manifolds, $q\colon V^n \to M^n$, which is a
homeomorphism away from $q^{-1}(\{p\})$, and for which
$q^{-1}(\{p\})$ is a nonempty compact set, classically a
projective space.  Real and complex versions of the construction 
are both considered here.
Some other notions of blowing up points appear in dynamics,
notably in studies of normal forms for vector fields and in 
constructions which delete a fixed point set and manipulate an
open cylinder.  

The principal theme presented here is that the algebraic geometers'
form of blowup (Section \ref{Section:Model}) is so natural that it
easily induces continuous homomorphisms on diffeomorphism groups
(Sections \ref{Section:Natural} and \ref{Sec:Maps}), defined by an
explicit model given in Section \ref{Section:Natural}. These lifting
homomorphisms render derivative data at the space level since the
exceptional locus in blowup is a projectivized tangent space.  For
complex manifolds and biholomorphic maps this rendering works nicely,
as it does for real $C^\infty$ manifolds and diffeomorphisms.
Dynamical consequences of the construction are laid out in Section
\ref{Sec:Dynamics}.

However, for finitely differentiable diffeomorphisms the
loss of regularity (Example \ref{Example:Regularity})
in blowup
becomes interesting and leads to a second theme: at
the lowest order of differentiability we find that a $C^1$
diffeomorphism fixing a point might lift to many homeomorphisms with
variant dynamics and quotient projections (Section \ref{Section:Variant}).  
The argument for this
nonuniqueness claim uses local $C^0$ conjugacy facts for hyperbolic
fixed points of diffeomorphisms and suggests that the dynamical
universality of the classical blowup is much more distinctive than the
spatial or single--map aspects of the construction.  This $C^0$
variation of blowup also indicates limitations on neighborhood--based
invariants of dynamics.  
The paper's third theme is that while the blowup notion is easy
to generalize in a $C^0$ context (Section \ref{Section:Blowups}),
greatly enlarging on the topological effects of classical
blowup (Section \ref{Sec:Topology}),
when we give up differentiability entirely
it turns out that the only reasonable generalized blowups which allow
every homeomorphism of the base manifold to lift are necessarily
homeomorphic to the base manifold (Theorem \ref{Thm:NoTopLift}).

Topologists are familiar with blowup as a construction tool and
stabilizing device in four--manifold topology.  Nash \cite{Nash52}
posed questions, since amplified by others, on the equivalence
relation on manifolds which blowup generates.  Major progress on
Nash's space--level question was made in
\cite{AbkulutKing91,BenedettiMarin92}, and especially in the work of
Mikhalkin \cite{Mikhalkin97}.  Blowup equivalence of diffeomorphisms
or group actions may be ripe for study after those advances.

We close the introduction with some notational conventions.
$B^n$ denotes the open ball $\{{\mathbf x} \in \Reals^n : 
|{\mathbf x}| < 1\}$, while the closed disk
$D^n = \{{\mathbf x} \in \Reals^n : 
|{\mathbf x}| \leq 1\}$.  
The derivative of $h$ at $p$ is written $Dh|_p$.
The projection from a Cartesian product onto its $j$--th factor
is denoted $\projection_j$.

Mapping spaces for pairs appear frequently below.  If $M$ is a
manifold and $A \subset M$ then $\operatorname{Aut}_{C^{k}}(M,\ A)$
denotes the space of $C^k$ diffeomorphisms $f$ of $M$ such that $f(A)
= A$; $f$ is not obliged to fix $A$ pointwise.  The analogous reading
is used for spaces of homeomorphisms or other self--maps such as
$\operatorname{Homeo}(M, A)$.

\myendinput


\section{Blowups of a Manifold or Map} \label{Section:Blowups}

This section considers nonclassical, merely continuous versions of 
the notion of blowing up a manifold or a map between manifolds.
Although these are easy to construct topologically, they do not
ordinarily have the universality properties of the classical
constructions, and it seems likely that homomorphisms such as
those exhibited in Theorem \ref{Theorem:Homom} are a distinguishing
feature of the classical constructions whose description
begins  in  Section \ref{Section:Model}.

\begin{definition} \label{Def:Blowup}
A \emph{topological blowup} of an $n$--manifold $M^n$ at a point 
$p \in M$ is a quotient map $q\colon V^n \to M^n$ such that
\newline
(1) $V^n$ is also an $n$--manifold,
\newline
(2) $\Sigma := q^{-1}(\{p\})$ is a connected, compact,
nonempty subset of $V$, and
\newline
(3) $q|_{V \mysetminus \Sigma}\colon V \mysetminus \Sigma
\to M \mysetminus \{p\}$ is a homeomorphism.

A topological blowup of a self--map $f\colon M^n \to M^n$
of an $n$--manifold $M^n$ at a fixed point $p = f(p) \in M$ 
of $f$ is a topological blowup $q\colon V^n \to M^n$ 
of $M$ at $p$ together with a self--map ${\widetilde f}\colon V \to V$ 
such that this diagram commutes:
\begin{equation*}
\begin{CD}
V    @> {\widetilde f} >>   V \\
@VV q V      @VV q V \\
M  @>> f >     M,
\end{CD}
\end{equation*}
i.e., $q\circ {\widetilde f} = f \circ q \colon V \to M$.
\end{definition}

$\Sigma = q^{-1}(\{p\})$ is called the \emph{exceptional locus} and
$M$ is sometimes described as a ``blowdown'' of $V$.  Our first
examples are constructed top--downwards, by beginning with $V$
and $\Sigma$.

One could very reasonably add to Definition \ref{Def:Blowup} the
requirement that $V \mysetminus \Sigma$ should be dense in $V$.  We
shall not do so in this paper, but note in advance the relevance of
this density condition in Theorem \ref{Thm:NoTopLift}.

A subset $\Sigma$ of a manifold $V^n$ is
\emph{cellular} if there are closed sets $S_i \subset V^n$
such that $S_1 \supseteq S_2 \supseteq \dots \supseteq S_i
\supseteq S_{i+1} \dots$, 
$\Sigma = \cap_1^\infty S_i$ and
for every $i$, $S_i \cong D^n$ is a disk imbedded with bicollared
boundary.

\begin{example} \label{Example:CE}
If $\Sigma$ is a cellular subset of $V^n$ then 
$M = V/\Sigma$ is a manifold and the quotient map 
$q\colon V \to V/\Sigma$ is a blowup of $M$ at 
the image of $\Sigma$.
\end{example}

\begin{figure}[tb]
\epsfbox{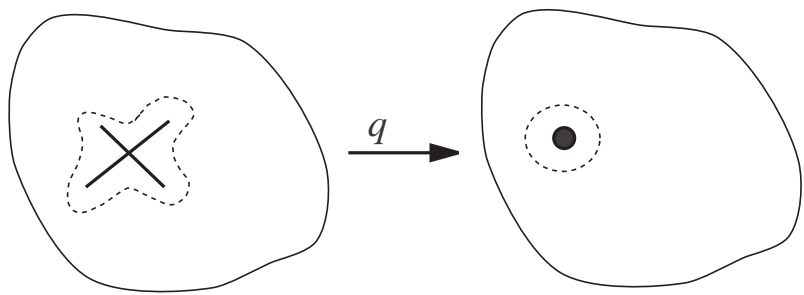}
\caption{Blowup with a cellular exceptional locus.}
\label{Figure:CE}
\end{figure}

See Figure \ref{Figure:CE} for a sketch of this sort of example, in
which $\Sigma$ is a finite polyhedral tree.

\epsfysize1.7in

\begin{example} \label{Example:MapCE}
If $g\colon V^n \to V^n$ is a map which preserves a cellular
subset $\Sigma \subset V$, then $g$ descends to a map
$f\colon V/\Sigma \to V/\Sigma$ and $g$ together with
the quotient map
$q\colon V \to V/\Sigma$ defines a blowup of $f$.
\end{example}

Examples \ref{Example:CE} and \ref{Example:MapCE} are misleading,
inasmuch as $V^n$ and $M^n$ are homeomorphic.  This is not usually
the case, and the replacement of $\{p\}$ by $\Sigma$ can affect
the global topology of a manifold in drastic ways.

\begin{example} \label{Example:Skeleton}
If $V^n \cong W^n \# X^n$ is a connected sum and
$\Sigma = X^{(n-1)}$ is the codimension--one skeleton
of a CW structure for $X$ which has one top--dimensional
cell, then the quotient map
$q\colon V \to V/\Sigma \cong W$ is a blowup of $W$
with exceptional locus $\Sigma$.
\end{example}

For instance, if $V^n$ is a compact, connected manifold and 
$\Sigma$ is the codimension--one skeleton of a cell structure for
$V$ which has one top cell then $V/\Sigma \cong S^n$.

\begin{figure}[tb]
\epsfbox{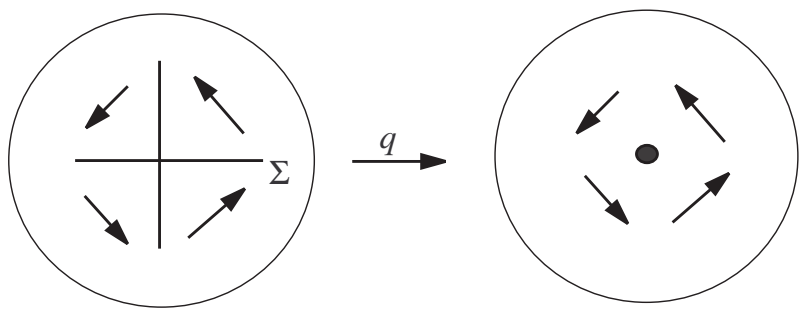}
\caption{Blowing up a rotation.}
\label{Figure:Rotation}
\end{figure}

\begin{example}
Examples of topological blowups for maps as well as spaces are not
hard to produce, and one is sketched in Figure
\ref{Figure:Rotation}.  Suppose that $g\colon V^n \to V^n$ is
periodic of period $r$ (so $g^r = \operatorname{Id}_V$), 
that $g$ has a fixed point $a$, and that the action of the
cyclic group $C_r$ generated by $g$ is effective on $V$ and
locally linear at $a$.  Form a
$g$--invariant polyhedral tree $\Sigma$ with $r$ legs emanating from
$a$, beginning with
a short segment $J$ based at $a$ so that $J \mysetminus \{a\}$
lies in the open dense set of $V$ on which $C_r$ acts freely.  
If $J$ is sufficiently short and becomes smooth in a linear model
for the action near $a$, then $\Sigma = \cup_1^r g^i(J)$ is the
desired tree, lying in a Euclidean ball about $a$.
The periodic map $g$ descends to a periodic map $f$ on $V/\Sigma$
and the pair $q\colon V \to V/\Sigma$, $g\colon V \to V$
defines a blowup of $f$.
\end{example}

\epsfysize1.7in

\begin{remark} \label{Remark:Collar}
The distinctive property of the exceptional locus $\Sigma$ in a 
blowup $q\colon V^n \to M^n$ arises from the requirement that
$V^n \mysetminus \Sigma \cong M^n \mysetminus \{p\}$ and
concerns deleted neighborhoods:  There is an
open neighborhood $U$ of $\Sigma$ in $V$ such that
$U \mysetminus \Sigma \cong B^n \mysetminus \{{\mathbf 0}\}$.

Such a neighborhood provides collared codimension--one spheres
exhibiting a connected sum structure for $V$, so Example
\ref{Example:Skeleton} is more typical than it might appear, although
up to this point we have allowed non--CW compacta to appear as
exceptional loci.  The classical blowup construction of algebraic
geometry exploits an instance of this neighborhood structure in
projective space, exactly along the lines of Example
\ref{Example:Skeleton}.
\end{remark}

\myendinput

\section{The Classical Model Construction} \label{Section:Model}
The most classical form of blowing up is performed at the origin in
$\Field^n$, where $\Field$ is $\Reals$ or $\Complexes$.  Our account
mostly follows \cite{McDuffSalamon98}, and a good description of the
construction and the properties which extend it from the affine model
to other varieties is found in \cite{Harris92}.

$\Projective(\Field^n)$ denotes the projective space of the
vector space $\Field^n$, defined as the quotient 
$\Projective(\Field^n) = (\Field^n \mysetminus \{{\mathbf 0}\})/\sim$,
where ${\mathbf v} \sim {\mathbf w}$ if and only if there exists
$\lambda \in \Field \mysetminus \{0\}$ such that 
${\mathbf v} = \lambda{\mathbf w}$.  Square brackets denote homogeneous
coordinates on a projective space, so that
the image in $\Projective(\Field^n)$
of $(v_1, \dots, v_n) \in \Field^n \mysetminus \{{\mathbf 0}\}$
is written $[v_1,\dots,v_n]$.
We will also use $[{\mathbf v}]$ to label the image in
$\Projective(\Field^n)$ of a nonzero vector $\mathbf v$ in
$\Field^n$.

Let $X \subset \Field^n \times \Projective(\Field^n)$ be the
subset
\begin{equation*}
X = \{\left((x_1,\dots,x_n), [y_1,\dots,y_n]\right) : 
     \ \text{for every $j,k$,}
     \ x_j y_k = x_k y_j \}
\end{equation*}
and let 
\begin{align*}
q\colon X &\to \Field^n \\
({\mathbf x}, [{\mathbf y}]) &\mapsto {\mathbf x}
\end{align*}
be the restriction of first--coordinate projection 
$\operatorname{pr}_1\colon \Field^n \times \Projective(\Field^n)
\to \Field^n$.

\begin{lemma}
$X =  \{({\mathbf x}, [{\mathbf y}]) : \ \text{there exists}
\ \mu \in \Field\ \text{such that}\ {\mathbf x} = \mu{\mathbf y} \}$.
In addition, $X$ is
a subvariety of $\Field^n \times \Projective(\Field^n)$,
$q$ is an algebraic map,  and preimages under $q$,
\begin{equation*}
q^{-1}(\{{\mathbf x}\}) =
\begin{cases}
\{ ({\mathbf x}, [{\mathbf x}]) \},& 
   \text{if ${\mathbf x} \neq {\mathbf 0}$},\\
\{ ({\mathbf 0}, [{\mathbf y}]) : [{\mathbf y}] \in 
   \Projective(\Field^n)\},&
   \text{ if ${\mathbf x} = {\mathbf 0}$},
\end{cases}
\end{equation*}
are such that $q$ is an isomorphism away from the origin 
and the fiber of $q$ over the origin is isomorphic to
the projective space $\Projective(\Field^n)$.
\qed
\end{lemma}

$\Sigma = q^{-1}(\{{\mathbf 0}\}) \cong \Projective(\Field^n)$ 
is usually called the
\emph{exceptional locus} or \emph{exceptional divisor}.
The quotient map $q$ is sometimes called the \emph{blowdown}
map, since it alters $X$ only by identifying $\Sigma$ to a point
(thus ``blowing down $\Sigma$'').

First--coordinate projection in $\Field^n \times
\Projective(\Field^n)$ defined the blowdown map $q$, and
second--coordinate projection determines the structure of a
neighborhood of the exceptional locus in classical blowups.

\begin{lemma} \label{Lemma:LineBundle}
Second--coordinate projection restricts
to $X$ as 
\begin{align*}
({\mathbf x}, [{\mathbf x}]) &\mapsto
[{\mathbf x}] \\
({\mathbf 0}, [{\mathbf y}]) &\mapsto
[{\mathbf y}]
\end{align*}
and identifies $X$ with the universal 
line bundle over $\Projective(\Field^n)$, i.e., the $\Field^1$--bundle
over this projective space whose fiber at $[{\mathbf y}]$ is the
line $\{\lambda{\mathbf y} : \lambda \in \Field\}$ 
through the origin and $\mathbf y$ in $\Field^n$.  
$\Sigma$ is identified with the zero section in this bundle, so
the normal bundle of $\Sigma$ in $X$ is identified with
the universal line bundle.
\qed
\end{lemma}

We will return to this bundle structure in Section \ref{Sec:Topology}.

\myendinput

\section{Naturality Properties and Manifold Constructions}
\label{Section:Natural}
Any self--map of $(\Field^n, \{{\mathbf 0}\})$
such that $h$ is differentiable at the origin and
$Dh|_0$ is a $\Field$--linear isomorphism
lifts to a map ${\widehat h}$ of $X$,
where ${\widehat h}$ is defined by:
\begin{align}  \label{Equation:HatMapDef}
{\widehat h}\colon ({\mathbf x}, [{\mathbf x}]) &\mapsto
(h({\mathbf x}), [h({\mathbf x})]), \\
\notag
{\widehat h}\colon ({\mathbf 0}, [{\mathbf y}]) &\mapsto
({\mathbf 0}, [Dh|_{\mathbf 0}({\mathbf y})]).
\end{align}

\begin{lemma}
If $h\colon (\Field^n, \{{\mathbf 0}\}) \to
(\Field^n, \{{\mathbf 0}\})$ is a continuous map which is
differentiable at ${\mathbf 0}$, and 
if $Dh|_0$ is a $\Field$--linear isomorphism, then
the map ${\widehat h}\colon (X, \Sigma) \to (X, \Sigma)$
defined above is continuous and makes this diagram commute
\begin{equation*}
\begin{CD} (X, \Sigma) @> {\widehat h} >> (X, \Sigma) \\
@V q VV                      @VV q V \\
(\Field^n, \{{\mathbf 0}\}) @>h >> (\Field^n, \{{\mathbf 0}\}),
\end{CD}
\end{equation*}
i.e., $q \circ {\widehat h} = h \circ q\colon
(X, \Sigma) \to (\Field^n, \{{\mathbf 0}\})$.
\end{lemma}

\begin{proof}
The claim that
$q \circ {\widehat h} = h \circ q$ follows immediately from
$q = \operatorname{pr}_1|_X$ and ${\widehat h}(\Sigma) = \Sigma$.
Because $q\colon X \mysetminus \Sigma \to \Field^n \mysetminus 
\{{\mathbf 0}\}$ is a homeomorphism and $h$ is continuous, the
continuity claim only needs to be confirmed at points of $\Sigma$.

The restriction of $\widehat h$ to $\Sigma = \{{\mathbf 0}\} \times
\Projective(\Field^n)$ is the projectivization
$\Projective(Dh|_{\mathbf 0})$ of a linear isomorphism, so this
restriction is a $C^\infty$ diffeomorphism on $\Sigma$.  Because
$h({\mathbf x}) = Dh|_{\mathbf 0}({\mathbf x}) + R({\mathbf x})$,
where $R({\mathbf x}) = o(|{\mathbf x}|)$ as
$|{\mathbf x}| \to 0$, $\widehat h$ is continuous
on the normal line $\{(\lambda {\mathbf x}, [\lambda {\mathbf x}])\}$
through $({\mathbf 0}, [{\mathbf x}]) \in \Sigma$:
\begin{align*}
{\widehat h}\left(\lambda {\mathbf x}, [\lambda {\mathbf x}]\right) &= 
\left(h(\lambda {\mathbf x}), [h(\lambda {\mathbf x})]\right) \\ &= 
\left(h(\lambda {\mathbf x}), [Dh|_{\mathbf 0}(\lambda {\mathbf x})
+ R(\lambda {\mathbf x})]\right) \\ &= 
\left(h(\lambda {\mathbf x}), 
[(\lambda|{\mathbf x}|)^{-1}Dh|_{\mathbf 0}(\lambda {\mathbf x})
+ (\lambda|{\mathbf x}|)^{-1} R(\lambda {\mathbf x})]\right) \\ &= 
\left(h(\lambda {\mathbf x}), 
[Dh|_{\mathbf 0}(|{\mathbf x}|^{-1}{\mathbf x})
+ (\lambda|{\mathbf x}|)^{-1} R(\lambda {\mathbf x})]\right),
\end{align*}
which tends to $({\mathbf 0}, [Dh|_{\mathbf 0}({\mathbf x})])$
as $\lambda \to 0$, with convergence uniform 
in $|{\mathbf x}|^{-1}{\mathbf x}$.
Therefore, by the triangle
inequality, $\widehat h$ is continuous at every point of $\Sigma$.
\end{proof}

The next lemma follows from the Chain Rule.  Recall that a map or
homeomorphism of pairs $(X, \Sigma) \to (X, \Sigma)$ is required
to carry $\Sigma$ to itself but need not restrict to the identity
on $\Sigma$.

\begin{lemma}  \label{Lemma:Homomorphism}
Let $g,h\colon (\Field^n, \{{\mathbf 0}\}) \to
(\Field^n, \{{\mathbf 0}\})$ be continuous maps which are 
differentiable at ${\mathbf 0}$ and have $\Field$--linear 
isomorphisms as their derivatives at the origin.  Then
$\widehat{g \circ h} = {\widehat g} \circ {\widehat h}\colon
(X, \Sigma) \to (X, \Sigma)$.
\qed
\end{lemma}

Since ${\widehat {\Identity_{\Field^n}}} = \Identity_X$,
the Lemma shows that $h \mapsto \widehat h$ defines a homomorphism
\begin{equation*}
\beta\colon\operatorname{Aut}_{C^{1}, \Field}(\Field^{n},\ \{{\mathbf 0}\}) \to
\operatorname{Homeo}(X,\ \Sigma),
\end{equation*}
where the $\Field$ decoration indicates that the derivatives are
required to be $\Field$--linear.

\begin{prop}
A smooth real or complex manifold $M$ can be blown up at any point $p$
to produce a quotient map from a smooth real or complex manifold
$\widehat M$,
\begin{equation*}
q\colon ({\widehat M}, \Sigma) \to (M, \{p\}),
\end{equation*}
which restricts to an isomorphism
$q|\colon {\widehat M} \mysetminus \Sigma \xrightarrow{\cong}
M \mysetminus \{p\}$.
If $M$ is modelled on $\Field^n$ then 
$\Sigma \cong \Projective(\Field^n)$.
\end{prop}

\begin{proof}
Lemma \ref{Lemma:Homomorphism} shows that origin--preserving
coordinate changes
with $\Field$--linear derivatives act as automorphisms of
$q\colon (X, \Sigma) \to (\Field^n, \{{\mathbf 0}\})$.
If $\phi_i\colon U \to \Field^n$ are local coordinate systems
($i=1,2$)
on a neighborhood $U$ of $p$ in $M$ then 
${\widehat {(\phi_1\circ\phi_2^{-1})}} \colon X \to X$
gives a change of coordinates on the model blowup.  The 
homomorphism properties established in
Lemma \ref{Lemma:Homomorphism} show that blownup coordinate change
maps satisfy the cocycle condition, yielding a consistent pasting
construction for $\widehat M$ from the data defining $M$.
\end{proof}

The same naturality properties used above give diffeomorphism blowups
on smooth manifolds, which are treated in detail in
Theorem \ref{Theorem:Homom}.

\myendinput


\section{Topology} \label{Sec:Topology}

This section describes the effects of the classical blowup
construction on topology.  We begin with the model construction at the
origin in $\Field^n$, where the space $X$ and the normal bundle of
$\Sigma$ in $X$ are identified with the universal line bundle over the
projective space $\Sigma$.

\begin{figure}[tb]
\epsfbox{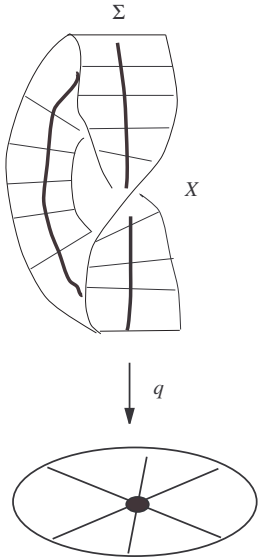}
\caption{Normal structure of the exceptional locus.}
\label{Figure:Mobius}
\end{figure}

This bundle description gives a picture of the blowdown map which may
be helpful (see Figure \ref{Figure:Mobius}).  Let $S(\Field^n)$ denote
the unit sphere in $\Field^n$; then a tubular neighborhood of $\Sigma$
in $X$ is identified with the mapping cylinder of the Hopf map
$h\colon S(\Field^n) \to \Projective(\Field^n)$ and the blowdown
quotient on this tubular neighborhood is the natural map between the
mapping cylinders for this Hopf map and for the constant map
$c\colon S(\Field^n) \to \text{point}$, i.e., the map of
pairs $(\operatorname{MapCyl}(h), \Projective(\Field^n)) \to
(\operatorname{MapCyl}(c), \{\text{point}\})$.

\epsfysize3in

In the complex case a bit of attention is required to the line bundles
playing roles in this discussion.  McDuff and Salamon
\cite{McDuffSalamon98} describe these identifications or computations
carefully:
\newline
(a) $\nu_X(\Sigma)$ is identified with the universal
line bundle $L$ over the projective space 
$\Sigma \cong \Projective(\Complexes^n)$;
\newline
(b) the first Chern class $c_1(\nu_X(\Sigma) = -c$, where
$c$ is the positive or canonical generator of 
$H^2(\Sigma; \Integers)$;
\newline
(c) the  normal line bundle to the hyperplane section in
$\Projective(\Complexes^{n+1})$ has
first Chern class $c_1(\nu_{\Projective(\Complexes^{n+1})}
\Projective(\Complexes^n)) = c$; and
\newline
(d) the  normal line bundle to the hyperplane section in
the conjugate complex structure
$\overline{\Projective(\Complexes^{n+1})}$ has
first Chern class 
$c_1\left(\nu_{\,\overline{\Projective(\Complexes^{n+1})}\,}
\overline{\Projective(\Complexes^n)}\right) = -c$.

The real blowup of a point in a Riemann surface has
$\Sigma \cong \Reals P^1 \cong S^1$, where the model space
$X$ is the nonorientable line bundle over $S^1$ whose total space
is a M\"obius band.  Thus,
for surfaces the mapping cylinder description of blowing up
and down suggests that blowing up a point has the global topological
effect of sewing in a crosscap.  This is true, and in general
the global effect of blowing up a point is a connected sum operation,
as in Example \ref{Example:Skeleton} and Remark \ref{Remark:Collar}.
For real blowups,
\begin{equation*}
{\widehat M}^n \cong M^n \# \Reals{}P^n,
\end{equation*}
and for complex blowups
\begin{equation*}
{\widehat M}^n \cong M^n \# \overline{\Complexes{}P^n}.
\end{equation*}
A conjugate complex structure appears in the second connected
sum because of the determinations of line bundles in the preceding
paragraph.

\myendinput

\section{Blowing Up Maps at a Fixed Point} \label{Sec:Maps}

The naturality properties of the classical blowup construction
suggest that (\ref{Equation:HatMapDef}) shows how to extend blowup
to a homomorphism of diffeomorphism groups.  This is possible,
but a kink develops in the $C^k$ case.

The regularity loss in the theorem below is formally due to a division
when one considers the homogeneous coordinate side of the formula for
$\widehat h$.   More geometrically, the blowup construction renders
tangential data for $h$ as spatial data for
$\widehat h$ since $\Sigma$ is the space of
lines in $T_pM$, so the loss of one derivative should be expected.
An example is worked out below to show that the
the loss of a derivative is genuine.

\begin{theorem}                    \label{Theorem:Homom}
Classical blowup of a point $p \in M$ determines
continuous, injective homomorphisms
\begin{equation*}
\beta\colon\operatorname{Aut}_{C^{k}}(M^{n},\ \{p\}) \to
\operatorname{Aut}_{C^{k-1}}(\widehat{M},\ \Sigma)
\end{equation*}
(in the real case) and
\begin{equation*}
\beta\colon\operatorname{Aut}_{\text{Holo}}(M^{n},\ \{p\}) \to
\operatorname{Aut}_{\text{Holo}}(\widehat{M},\ \Sigma)
\end{equation*}
(in the complex analytic case), where
$\beta(h) = {\widehat h}$ in both cases.
\end{theorem}

\begin{proof}
Away from $\Sigma$ we know that 
$\widehat h$ and $h$ may be identified, so the regularity issue
only arises along $\Sigma$.  Along $\Sigma$ we have defined
$\widehat h$ to be the projectivization of the linear
map $Dh|_{\mathbf 0}$, so $\widehat h$ is infinitely differentiable
in those directions.

If for ${\mathbf x}$ near ${\mathbf 0}$ we have
$h({\mathbf x}) 
= \sum_{j=1}^k D^j h|_{\mathbf 0} ({\mathbf x},\dots, {\mathbf x})
+ R({\mathbf x})$, where $R({\mathbf x}) = o(|{\mathbf x}|^k)$
as ${\mathbf x} \to {\mathbf 0}$, then
\begin{align*}
{\widehat h}({\mathbf x}, [{\mathbf x}]) 
&= 
\left( \sum_{j=1}^k D^j h|_{\mathbf 0} ({\mathbf x},\dots, {\mathbf x})
+ R({\mathbf x}),\ %
[\sum_{j=1}^k D^j h|_{\mathbf 0} ({\mathbf x},\dots, {\mathbf x})
+ R({\mathbf x})] \right)\\
&= 
\left( \sum_{j=1}^k D^j h|_{\mathbf 0} ({\mathbf x},\dots, {\mathbf x})
+ R({\mathbf x}),\ %
[\sum_{j=1}^k D^j h|_{\mathbf 0} ({\mathbf x/|{\mathbf x}|},{\mathbf x},
\dots, {\mathbf x})
+ R({\mathbf x})/|{\mathbf x}|]
\right).
\end{align*}
The division by $|{\mathbf x}|$ inside the homogeneous coordinates
gives a zero--th order term of 
$Dh|_{\mathbf 0}({\mathbf x}/|{\mathbf x}|)$, similarly reduces the
degree of the other homogeneous terms in the Taylor expansion,
and reduces by one the order of vanishing for the remainder term.
The resulting expansion of $\widehat h$ near 
$(0,\ [{\mathbf x}]) \in \Sigma$ shows that we lose one
partial derivative
of $\widehat h$ along the fibers of the normal bundle $\nu_X(\Sigma)$,
compared to the degree of smoothness of $h$ at ${\mathbf 0}$.
(Recall from Section \ref{Section:Natural}
 that $\widehat h$ is $C^\infty$ along $\Sigma$.)

The partial derivatives of $\widehat h$ along the singular locus
and normal to it are continuous, through order $k-1$, so $\widehat h$
is $C^{k-1}$ at points of $\Sigma$ by the familiar theorem deducing
(Fr\'echet) differentiability from continuous partial derivatives.

Lemma \ref{Lemma:Homomorphism} and surrounding discussion show that 
$\beta$ is a homomorphism.

Once we know that $\widehat h$ is continuous, it
follows that  $\beta$ is injective,
since $h$ determines $\widehat h$ on the dense subset
${\widehat M} \mysetminus \Sigma$.

$\beta$ is continuous because the $C^r$ distance between
diffeomorphisms on $M^n$ majorizes the $C^{r-1}$ distance between
their blowups on ${\widehat M} \mysetminus \Sigma$ and the $C^{r-1}$
distance between those blowups on $\Sigma$.
\end{proof}

\begin{example}          \label{Example:Regularity}
This is a two--dimensional example of
the regularity loss from $C^1$ to $C^0$
indicated in the theorem.

Let $g\colon \Reals \to \Reals$ be defined by
$g(x) = x + x|x|$, so 
that $g'(x) = 1 + 2|x|$.  $g$ is a $C^1$ diffeomorphism, 
but not $C^2$, and $g(0) = 0$.

Define $h\colon \Reals^2 \to \Reals^2$ to be the $C^1$ diffeomorphism
$h(x,y) =  (g(x), y)$.  This map preserves the origin and blows up there
to ${\widehat h}\colon \left( (x,y),\ [x,y] \right) \mapsto
\left( (x + x|x|, y),\ [x+ x|x|, y] \right)$.  The parametrized
line
$t \mapsto (t, mt)$ of slope $m \neq 0$ in the plane is covered in
the blownup plane by the $C^{\infty}$ parametric
curve $c\colon t \mapsto \left((t, mt),\ [t, mt] \right)$, 
and the composite
\begin{align*}
{\widehat h} \circ c\colon t 
&\mapsto
\left( (t + t|t|, mt),\ [t+ t|t|, mt] \right) \\
&= 
\left( (t + t|t|, mt),\ [m^{-1} + m^{-1}|t|, 1] \right) \\
\end{align*}
is continuous but not differentiable at $t = 0$, since in the
usual local coordinate system about $[0,1] \in \Projective(\Reals^2)$
the second component becomes $t \mapsto m^{-1} + m^{-1}|t|$:
therefore $\widehat h$ is not differentiable, and 
not even G\^ateaux differentiable,
at $c(0) = \left((0,0),\ [1,m]\right) \in \Sigma$.
\end{example}

Similar examples for $C^k$ to $C^{k-1}$ regularity loss are
available for all $k \geq 1$.

Theorem \ref{Thm:NoTopLift}
indicates that $C^0$ to $C^0$ lifting of automorphisms through
blowups is problematic for other
reasons.  

\myendinput

\section{Dynamics} \label{Sec:Dynamics}

The dynamics of a classically blown--up diffeomorphism $\widehat h$
off, on, and near the singular locus are described in terms
of basic features of the original diffeomorphism $h$.

If $h\colon (M, \{p\}) \to (M, \{p\})$ is a diffeomorphism
fixing $p$ then the restriction of $q$ gives
\begin{equation*}
{\widehat h}|_{{\widehat M} \smallsetminus \Sigma}
\cong
h|_{ M \smallsetminus \{p\}},
\end{equation*}
so blowup does not modify dynamics far from the exceptional locus.

\begin{lemma}  \label{Lemma:BlowFix}
On the exceptional locus $\Sigma = \Projective(T_pM)$
\begin{equation*}
{\widehat h}|_{\Sigma} \cong
\Projective(Dh|_p)
\end{equation*}
is a projectivized linear map, with the new fixed point set 
given by
\begin{equation} \label{Equation:SigmaFixedSet}
\Sigma\,\cap\,\operatorname{Fix}({\widehat h}) = 
\amalg_{\lambda \in \Field}
\ \Projective(E_\lambda),
\end{equation}
where $E_\lambda = \ker(\lambda I - Dh|_p)$
and $\amalg$ denotes a disjoint union.
\end{lemma}

\begin{proof}
Equation \ref{Equation:HatMapDef} defines
${\widehat h}\colon ({\mathbf 0}, [{\mathbf y}]) \mapsto
({\mathbf 0}, [Dh|_{\mathbf 0}({\mathbf y})])$ in the
model case, so the restriction
of $\widehat h$ to $\Sigma$ is the projectivization of
the derivative $Dh|_p$.

Therefore, fixed points of $\widehat h$
in $\Sigma = \Projective(T_pM)$ are solutions of
$[Dh|_p(v)] = [v]$, i.e. projective
equivalence classes of tangent vectors $v$ for which there
exist scalars $\lambda \in \Field$ satisfying
$Dh|_p(v) = \lambda v$, that is, projective equivalence
classes of eigenvectors of $Dh|_p$.  
Equation \ref{Equation:SigmaFixedSet} describes the set of
all such projective classes as a disjoint union of projectivized
subspaces of $T_p M$.
\end{proof}

Derivative computations for $\widehat h$ at points of $\Sigma$ 
and in directions not tangent to $\Sigma$ will
involve the loss of order noted
in Theorem \ref{Theorem:Homom}.  These are omitted here -- see the
displayed equation in the proof of that theorem for the appearance of
second derivatives of $h$ in the first derivative of $\widehat h$.
Despite this derivative complication,
a qualitative picture of part of the dynamics of $\widehat h$ normal
to $\Sigma$ is provided by some naturality observations.  

First, if $N^k$ is a $C^1$ submanifold of $N^n$ passing through
$p$ and $\widehat M$ is the blowup of $M^n$ at $p$ then there is a
$C^0$ submanifold $V^k$ of $\widehat M$ such that $V$ is homeomorphic
to $\widehat N$ and $q_N \colon ({\widehat N}, \Sigma_N) \to
(N, \{p\})$ is equivalent to $q_M|\colon (V, V\cap \Sigma_M)
\to (N, \{p\})$.  The main step in checking this claim is a confirmation that 
the closure of $q_M^{-1}(N \mysetminus \{p\})$ in $\widehat M$ meets
$\Sigma_M$ in $\Projective(T_p N^k)$.

Second, if $h\colon (M^n, N^k, \{p\}) \to (M^n, N^k, \{p\})$ is
a $C^1$ diffeomorphism then $\widehat h$ preserves the submanifold
$V^k$ defined in the preceding paragraph.

\begin{figure}[tb]
\epsfbox{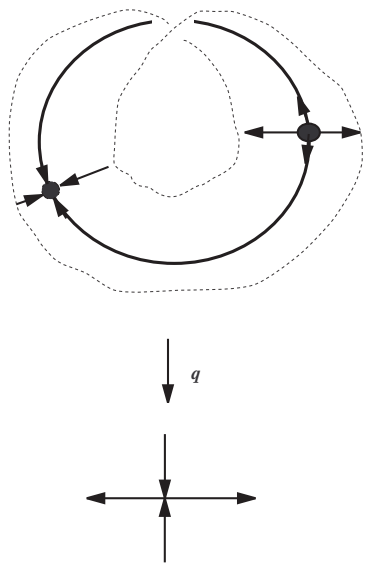}
\caption{Blowing up a hyperbolic fixed point.}
\label{Figure:Hyperbolic}
\end{figure}

In particular, if $p$ is a hyperbolic fixed point
of $h$ then this applies to the stable and unstable manifolds at $p$,
and also to any invariant local submanifolds tangent to 
other invariant subspaces of $T_p M$, such as eigenspaces of
$Dh|_p$.  Because $\Sigma \cap V = \Projective(T_pN)$, 
in Figure \ref{Figure:Hyperbolic} the blownup one--dimensional
stable submanifold at $p$ 
meets $\Sigma$ in the point corresponding to the appropriate
eigenspace of $Dh|_p$, while the blownup unstable manifold meets
$\Sigma$ in the point corresponding to another one--dimensional
eigenspace.

\epsfysize3in

This low--dimensional example suggests the behavior of blowups at
hyperbolic fixed points, but is a bit simpler than the general case,
which we sketch now.

Suppose that $p$ is a hyperbolic fixed point of $h$, that $E_s$ and
$E_u$ are the stable and unstable subspaces of $T_p M$, and that the
stable and unstable submanifolds at $p$ are $W_s$ and $W_u$.  $W_s$
and $W_u$ blow up at $p$ to give invariant submanifolds $V_s$, $V_u$
which meet $\Sigma$ in a submanifold (either $\Projective(E_s)$ or
$\Projective(E_u)$) which is invariant and forward attracting
(respectively backward attracting).
The dynamics of $\widehat h$ restricted to 
$\Projective(E_s)$ or
$\Projective(E_u)$) are those of a projectivized linear map.

\myendinput


\section{Variant Blowups} \label{Section:Variant}

$C^0$ data for $h$ are not enough to determine ${\widehat h}|_\Sigma$
as a homeomorphism.  For example, in the hyperbolic case we can apply
local conjugacy results to obtain lots of homeomorphisms blowing up a
given $C^1$ diffeomorphism.  A handy reference for these facts on
topological conjugacy is \cite[Sec.~6.3]{KatokHasselblatt95}, which
gives a proof of this result on local equivalence of hyperbolic fixed
points:

\begin{remark}  \label{Remark:Conjugacy}
Topological conjugacy classes of hyperbolic diffeomorphisms
with $p$ as an isolated fixed point are
determined by the dimensions and orientations of the stable
and unstable manifolds of these diffeomorphisms at $p$.
\end{remark}

\begin{figure}[tb]
\epsfbox{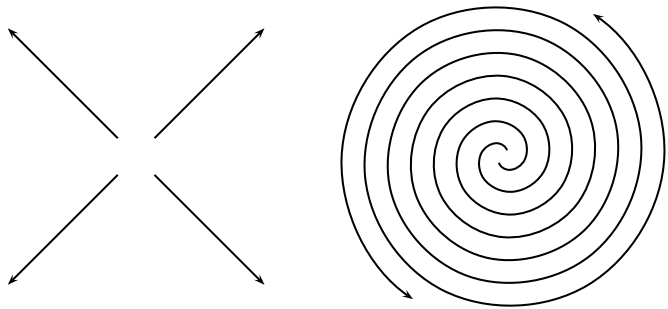}
\caption{Topologically conjugate hyperbolic fixed points.}
\label{Figure:Spiral}
\end{figure}

For example, if $M$ is even--dimensional over $\Reals$, $p$ is a fixed
point for $h_i$ ($i=1,2$), and $D{h_1}|_p$ is diagonalizable with all
eigenvalues lying in the interval $\lambda > 1$, while $D{h_2}|_p$ has
only non--real eigenvalues, all satisfying $|\lambda| > 1$, then $h_1$
and $h_2$ are locally conjugate near $p$.  Figure \ref{Figure:Spiral}
suggests how different in appearance such topologically conjugate
diffeomorphisms can be, and indicates that the cause of the phenomenon
is a familiar difficulty: we can unwind a spiral with a continuous
automorphism, but not with a smooth one.

\epsfysize2in

A global topological conjugacy $\phi\colon M \to M$ from
$h_1$ to $h_0$ leads to a variant blowup of $h_0$ with
${\widehat h}_1$ as the covering homeomorphism and
$\phi\circ q$ as blowdown map.

\begin{equation*}
\begin{CD}
{\widehat M} @> {\widehat h}_1 >> {\widehat M} \\
@V q VV                             @VV q V \\
M            @> h_1 >>              M \\
@V \phi VV                          @VV \phi V \\
M            @> h_0 >>              M
\end{CD}
\end{equation*}

The dynamics on $\Sigma$ of these conjugacy--induced blowups can
differ dramatically from those of the classical construction.  Fixed
point sets on $\Sigma$ may differ drastically in dimension as we run
over diffeomorphisms $h_1$ which are topologically conjugate to $h_0$,
ranging from empty to discrete to connected and high--dimensional.
We emphasize this variability with a proposition.

\begin{prop}
Let $h\colon M^n \to M^n$ be a diffeomorphism with an isolated
hyperbolic fixed point $p$.  If the stable and unstable
subspaces $E^s_p, E^u_p \subseteq T_pM$ are both even--dimensional
then there are conjugacy--induced topological blowups 
${\widetilde h}$ of $h$ such that 
$\Sigma \cap \operatorname{Fix}({\widetilde h})$
is of any of these sorts:
\newline
(a) empty,
\newline
(b) discrete, containing any even number of points between $2$ and $n$,
or
\newline
(c) positive--dimensional, with any dimension between $1$ and 
$-1 + \max(\dim(E^u_p),\,\dim(E^s_p))$.
\end{prop}

\begin{proof}
In each case a $C^0$
conjugacy as described in Remark \ref{Remark:Conjugacy}
between $h$ and another diffeomorphism $g$ with a
hyperbolic fixed point at $p$ yields the topological blowup.
Our job here is to allocate eigenvalues for $Dg|_p$ and apply
Lemma \ref{Lemma:BlowFix} to $\widehat g$.

If every eigenvalue $\lambda \in \Complexes \mysetminus \Reals$
then $Dg|_p$ has no real eigenvectors and the classical blowup's
$\Sigma \cap \operatorname{Fix}({\widehat g})$ is empty.

If $Dg|_p$ has $k$ distinct real eigenvalues and $n -  k$ 
complex eigenvalues appearing in conjugate pairs, then 
$k$ must be even but may otherwise assume any value between
$0$ and $n$.  In this case we see $k$ isolated fixed points
for $\widehat g$ on $\Sigma$.

Positive--dimensional fixed point sets arise from repeated real
eigenvalues for $Dg|_p$.  These may appear in combinations so that
the multiplicities of real eigenvalues form partitions of some
even numbers $0 \leq e_u \leq \dim(E_p^u)$,
$0 \leq e_s \leq \dim(E_p^s)$:
\begin{align*}
k_1 + k_2 + \dots + k_q &= e_u,\\
l_1 + l_2 + \dots + l_r &= e_s,\\
\end{align*}
where $k_i$, $l_j$ are the multiplicities of unstable, respectively
stable, real eigenvalues
of $Dg|_p$.  Each of these real eigenvalues produces a component
of the fixed point set $\Sigma \cap \operatorname{Fix}({\widetilde g})$
which is diffeomorphic to a
projective space:  If $\lambda$ is a real eigenvalue of
multiplicity $m$ then $\Projective(E_{\lambda}) \cong
\Projective(\Reals^{m}) \cong \Reals P^{m - 1}$.
The largest dimension arising in this way
$\max(\dim(\Projective(E^u_p)),\,\dim(\Projective(E^s_p)))$.
\end{proof}

Note that $\Sigma \cap \operatorname{Fix}({\widehat g})$ might have
components of different dimensions.  The complex case is similar, but
the fixed point set must be nonempty.

The next few results indicate that
topological blowups are necessarily limited in naturality.

\begin{lemma} \label{Lemma:TubeHomeo}
Let $N^k$ be a connected topological manifold without boundary.
If $\{x_i\}$ and  $\{y_i\}$ are sequences in $N \times (0,\infty)$
such that $\{\projection_2(x_i)\} \to \infty$,
$\{\projection_2(y_i)\} \to \infty$, and both of these sequences
in $\Reals$ are strictly increasing, then there is a homeomorphism
$g\colon N \times (0,\infty) \to N \times (0,\infty)$ such that
for every $i$, $g(x_i) = y_i$, and such that for some $\varepsilon > 0$
the restriction
$g|_{(0,\varepsilon]}$ is the identity.
\end{lemma}

\begin{proof}
This is a consequence of the following version of the homogeneity
of manifolds: for any $x,y \in N$ there
exists an isotopy from $\operatorname{Id}_N$ to a homeomorphism
which carries $x$ to $y$.

In more detail, $g$ can be built in segments which are pasted
together.
We may apply a homeomorphism of the form $\rho\colon (x,t) \mapsto
(x, \psi(t))$, where $\psi\colon (0,\infty) \to (0,\infty)$ is
a homeomorphism, to arrange that $\projection_2(\rho(x_i))
= \projection_2(y_i)$ for all $i$.  Continue the argument with
the sequence $\{\rho(x_i)\}$ replacing $\{x_i\}$.

Let $0 < s_1 < \projection_2 (x_1)$ and
let $\phi_t\colon N \to N$ be an isotopy over $s_1 \leq t \leq 
\projection_2(x_1)$ such that $\phi_{s_1} = \operatorname{Id}_N$
and $\phi_{\projection_2(x_1)}(\projection_1(x_1)) = \projection_1(y_1)$.
Define the first two pieces of 
$g\colon N \times (0,\infty) \to N \times  (0,\infty)$ by 
$g|_{N \times (0,s_1]} = \operatorname{Id}_{N \times (0,s_1]}$ and
$g|_{N \times [s_1, \projection_2(x_1)]}\colon
(x,t) \mapsto (\phi_t(x), t)$.  
Subsequent segments are defined by taking an isotopy
$\phi_t$ over $\projection(x_{i-1})
\leq t \leq \projection_2(x_{i})$ which starts with
the $\phi_{\projection(x_{i-1})}$ already selected
and ends at a homeomorphism which carries 
$\phi_{\projection(x_{i-1})}(x_i)$ to $y_i$.
\end{proof}

The topological blowups $q\colon (V, \Sigma) \to (M, \{p\})$
of greatest interest will share with
the classical construction the property that 
$V \mysetminus \Sigma$ is dense in $V$.  This line of argument
shows that in no such case
can we find a lifting construction for homeomorphisms.

\begin{theorem} \label{Thm:NoTopLift}
Let $q\colon (V^n,\Sigma) \to (M^n, \{p\})$ be a topological blowup
of the manifold $M^n$ at $p$ such that at least two points lie on
the frontier of $\Sigma$ in $V$.  If $n \geq 2$
then there is a homeomorphism $h\colon (M, \{p\}) \to (M, \{p\})$
which does not lift through $q$.
\end{theorem}

\begin{proof}
Suppose that $x, y \in \Sigma$ are distinct points and that
$\{x_i\}$, $\{y_i\}$ are sequences in $V \mysetminus \Sigma$
such that $x = \lim_{i \to \infty} x_i$ and $y = \lim_{i \to \infty} y_i$.
We may assume that both sequences lie in a neighborhood $U$ of $\Sigma$
which admits a homeomorphism 
$f\colon U \mysetminus \Sigma \xrightarrow{\cong}
S^{n-1} \times (0,\infty)$ and that the sequences
$\{\projection_2\circ f(x_i)\}$ and 
$\{\projection_2\circ f(y_i)\}$ are strictly increasing and converge
to $\infty$.

Form a third sequence $\{z_i\}$ such that each $z_{2j}$ is one of the
$y_i$, each $z_{2j+1}$ is one of the $x_i$, and the real sequence
$\{\projection_2\circ f(z_i)\}$ is strictly increasing and converges
to $\infty$.  $\{z_i\}$ is not convergent in $V$, but all three of
the sequences $\{q(x_i)\}$, $\{q(y_i)\}$, and $\{q(z_i)\}$ converge
to $p$ in $M$.

Since $n \geq 2$, Lemma \ref{Lemma:TubeHomeo} 
implies that there is a homeomorphism $h\colon (M, \{p\}) \to 
(M, \{p\})$ such that for every $i$, $h(q(x_i)) = q(z_i)$.  If $h$
is covered by $\widetilde h\colon V \to V$ then 
$\{{\widetilde h}(x_i)\}$ is divergent although $\{x_i\}$ converges
to $x$, so $\widetilde h$ is not continuous.
\end{proof}

\begin{cor}
Suppose that $n \geq 2$ and let
$q\colon (V^n,\Sigma) \to (M^n, \{p\})$ be a topological blowup
of the manifold $M^n$ at $p$ such that $V \mysetminus \Sigma$ is
dense in $V$.  If every homeomorphism of $(M, \{p\})$ is
covered by a homeomorphism of $(V, \Sigma)$, then 
$\Sigma = \{\text{pt.}\}$ and $q$ is a homeomorphism.
\qed
\end{cor}

\myendinput
\providecommand{\bysame}{\leavevmode\hbox to3em{\hrulefill}\thinspace}

\end{document}